\documentclass[12pt]{article}
\usepackage{amssymb, amsmath, amsfonts, tikz}
\usepackage[indent,margin=1cm]{caption}
\usepackage[colorlinks,linkcolor=blue,anchorcolor=blue,citecolor=blue]{hyperref}
\newtheorem{thm}{Theorem}[section]

\newtheorem{prob}[thm]{Open problem}

\newtheorem{lem}[thm]{Lemma}
\newtheorem{coro}[thm]{Corollary}

\def\pf{\noindent{\it Proof.} }
\def\qed{\nopagebreak\hfill{\rule{4pt}{7pt}}
\medbreak}

\numberwithin{equation}{section}

\def\qed{\nopagebreak\hfill{\rule{4pt}{7pt}}
\medbreak}

\newlength{\boxedparwidth}
  {\begin{center} \begin{tabular}{|@{\hspace{.315in}}c@{\hspace{.15in}}|}
                  \hline \\ \begin{minipage}[t]{\boxedparwidth}
                  \setlength{\parindent}{.25in}}%
  {\end{minipage} \\ \\ \hline \end{tabular} \end{center}}

\begin{document}

\begin{center}
{\large \bf The $e$-positivity of two families of $(claw, 2K_2)$-free graphs}
%chromatic symmetric functions of
% generalized pyramids and $2K_2$-free unit interval graphs
\end{center}
\begin{center}
Grace M.X. Li$^{1}$ and Arthur L.B. Yang$^{2}$\\[6pt]

$^{1,2}$Center for Combinatorics, LPMC\\
Nankai University, Tianjin 300071, P. R. China\\[6pt]

Email: $^{1}${\tt limengxing@mail.nankai.edu.cn}, $^{2}${\tt yang@nankai.edu.cn}
\end{center}

\noindent\textbf{Abstract.} Motivated by Stanley's conjecture about the $e$-positivity of claw-free incomparability graphs, Hamel and her collaborators studied the $e$-positivity of $(claw, H)$-free graphs, where $H$ is a four-vertex graph. In this paper we establish the $e$-positivity of generalized pyramid graphs and $2K_2$-free unit interval graphs, which are two important families of $(claw, 2K_2)$-free graphs. Hence we affirmatively solve one problem proposed by Hamel, Ho\`{a}ng and Tuero, and another problem considered by Foley, Ho\`{a}ng and  Merkel.

\noindent \emph{AMS Classification 2010:} 05E05, 05C15


\section{Introduction}

Given a finite simple graph $G$ with vertex set $V$ and edge set $E$, a proper coloring of $G$ is a function $\kappa$ from $V$ to $\mathbb{P} = \{1,2,\ldots\}$ such that $\kappa(u)\neq\kappa(v)$ whenever $uv\in E$. Stanley \cite{Sta-95} defined the chromatic symmetric function $X_G$ as
\begin{align}\label{def-xg}
X_G = \sum_{\kappa} \prod_{v \in V} x_{\kappa(v)},
\end{align}
where $\kappa$ ranges over all proper colorings of $G$. It is clear that $X_G$ is a homogeneous symmetric function of degree $n$ if the cardinality of $V$ is $n$. There have been many works focusing on the expansion of $X_G$ in terms of various bases of symmetric functions. A well known basis is composed of elementary symmetric functions which are indexed by integer partitions.
Recall that an integer partition of $n$ is a weakly decreasing sequence $\lambda=(\lambda_1,\lambda_2,\ldots,\lambda_k)$ of positive integers such that $\sum_{i=1}^k\lambda_i=n$, denoted by $\lambda\vdash n$. Sometimes we consider $\lambda$ as an infinite sequence by appending infinite $0$'s.
The elementary symmetric function $e_{\lambda}$ is defined as
\begin{align*}
e_{\lambda}=e_{\lambda_1}e_{\lambda_2}\cdots e_{\lambda_k},
\end{align*}
where
\begin{align*}
e_0=1 \mbox{ and }e_i=\sum_{1\leq j_1<j_2<\cdots<j_i} x_{j_1}x_{j_2}\cdots x_{j_i} \mbox{ for }i\geq 1.
\end{align*}
It is well known that the set $\{e_{\lambda} \mid \lambda\vdash n\}$ forms a basis of homogeneous symmetric functions of degree $n$.
 A celebrated conjecture of Stanley \cite{Sta-95} states that the chromatic symmetric function $X_G$ of a claw-free incomparability graph $G$ is $e$-positive, namely, $X_G$ can be written as a nonnegative linear combination of $\{e_{\lambda}\}$. If $X_G$ is $e$-positive, we also say that $G$ is $e$-positive for convenience. Many works have been done towards Stanley's conjecture, see for instance \cite{Cho-Hon-19,Dah-18,Dah-She-Wil-19,Har-Pre-18,Sha-Wac-16}. The main objective of this paper is to prove the $e$-positivity of two families of $(claw, 2K_2)$-free graphs.

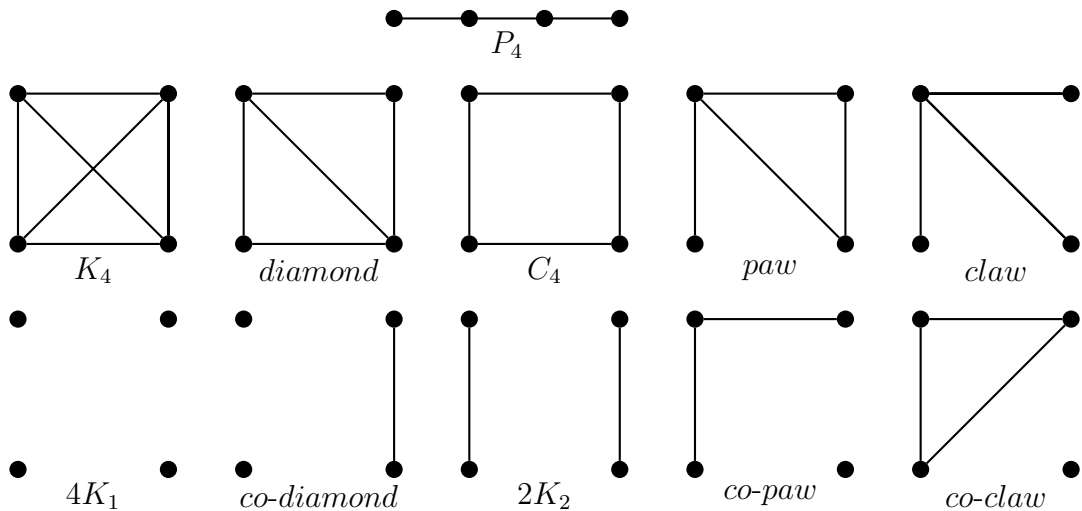
\begin{figure}[ht]
\center
\begin{tikzpicture}
[place/.style={thick,fill=black!100,circle,inner sep=0pt,minimum size=2mm,draw=black!100}]
\node [place] (v1) at (0,6) {};
\node [place] (v2) at (2,6) {};
\node [place] (v4) at (0,4) {};
\node [place] (v3) at (2,4) {};
\draw[thick] (v1) -- (v2) -- (v3) -- (v4) -- (v1) -- (v3) -- (v2) -- (v4);
\node at (1,3.65) {$K_4$};
\node [place] (v6) at (5,6) {};
\node [place] (v5) at (3,6) {};
\node [place] (v8) at (3,4) {};
\node [place] (v7) at (5,4) {};
\draw[thick] (v5) -- (v6) -- (v7) -- (v8) -- (v5) -- (v7);
\node at (4,3.65) {$diamond$};

\node [place] (v9) at (6,6) {};
\node [place] (v12) at (6,4) {};
\node [place] (v10) at (8,6) {};
\node [place] (v11) at (8,4) {};
\draw[thick] (v9) -- (v10) -- (v11) -- (v12) -- (v9);
\node at (7,3.65) {$C_4$};

\node [place] (v13) at (9,6) {};
\node [place] (v16) at (9,4) {};
\node [place] (v14) at (11,6) {};
\node [place] (v15) at (11,4) {};
\draw[thick] (v13) -- (v14) -- (v15) -- (v13) -- (v16);
\node at (10,3.65) {$paw$};

\node [place] (v17) at (12,6) {};
\node [place] (v20) at (12,4) {};
\node [place] (v18) at (14,6) {};
\node [place] (v19) at (14,4) {};
\draw [thick](v17) -- (v18) -- (v17) -- (v19) -- (v17) -- (v20);
\node at (13,3.65) {$claw$};

\node [place] at (0,3) {};
\node [place] at (0,1) {};
\node [place] at (2,3) {};
\node [place] at (2,1) {};
\node at (1,0.65) {$4K_1$};

\node [place] at (3,3) {};
\node [place] at (3,1) {};
\node [place] (v21) at (5,3) {};
\node [place] (v22) at (5,1) {};
\draw[thick] (v21) -- (v22);
\node at (4,0.65) {$co$-$diamond$};

\node [place] (v23) at (6,3) {};
\node [place] (v24) at (6,1) {};
\node [place] (v25) at (8,3) {};
\node [place] (v26) at (8,1) {};
\draw[thick] (v23) -- (v24);
\draw[thick](v25) -- (v26);
\node at (7,0.65) {$2K_2$};

\node [place] (v28) at (9,3) {};
\node [place] (v27) at (9,1) {};
\node [place] (v29) at (11,3) {};
\node [place] at (11,1) {};
\draw [thick](v27) -- (v28) -- (v29);
\node at (10,0.65) {$co$-$paw$};

\node [place] (v30) at (12,3) {};
\node [place] (v32) at (12,1) {};
\node [place] (v31) at (14,3) {};
\node [place] at (14,1) {};
\draw[thick] (v30) -- (v31) -- (v32) -- (v30);
\node at (13,0.65) {$co$-$claw$};

\node [place] at (5,7){};
\node [place] at (6,7){};
\node [place] at (7,7){};
\node [place] at (8,7){};
\draw[thick] (5,7) -- (6,7) -- (7,7) -- (8,7);
\node at (6.5,6.65) {$P_4$};
\end{tikzpicture}

\caption{All four-vertex graphs.}\label{gra-9}
\end{figure}

Let us first recall some related concepts and give an overview of some background. Let $H$ be a set of graphs. A graph $G$ is said to be $H$-free if it does not contain any graph of $H$ as an induced subgraph. Hamel, Ho\`{a}ng and Tuero \cite{Ham-Hoa-Tue-19} studied the $e$-positivity of $H$-free graphs,
where $H$ is composed of one claw and another four-vertex graph. There are eleven graphs on four vertices, see figure \ref{gra-9}.
Concerning the $e$-positivity of $(claw, F)$-free graphs with $F$ being a four-vertex graph other than claw, some progress has been made.
Tsujie \cite{Tsu-18} proved the $e$-positivity for the case $F=P_4$. Hamel, Ho\`{a}ng and Tuero proved the $e$-positivity for $F=paw$ and $F=co\mbox{-}paw$.
They also showed that a $(claw, F)$-free graph is not necessarily $e$-positive if $F$ is a diamond, co-claw, $K_4,\, 4K_1,\, 2K_2$ or $C_4$.
It remains to study the case $F$ is a co-diamond. Hamel, Ho\`{a}ng and Tuero proposed the following open problem.
\begin{prob}
Whether $(claw, co\text{-}diamond)$-free graphs are $e$-positive?
\end{prob}
By considering the structure of $(claw, co\text{-}diamond)$-free graphs, they reduced the above problem to determine the $e$-positivity of certain peculiar graphs, as illustrated in \cite[Figure 3]{Ham-Hoa-Tue-19}.

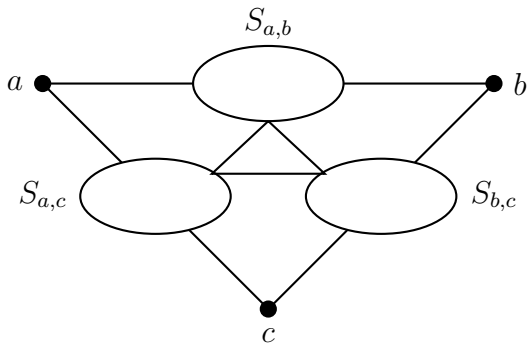
\begin{figure}[ht]
\center
\begin{tikzpicture}
[place/.style={thick,fill=black!100,circle,inner sep=0pt,minimum size=2mm,draw=black!100}]
\node [place, label=left:$a$] at (2,7) {};
\node [place, label=right:$b$] at (8,7) {};
\node [place, label=below:$c$] at (5,4) {};
\draw[thick] (2,7) -- (8,7) -- (5,4) -- (2,7);
\draw[thick][fill=white] (5,7) circle [x radius=1cm, y radius=5mm] node at (5,7.8){$S_{a,b}$};
\draw[thick][fill=white] (3.5,5.5) circle[x radius=1cm, y radius=5mm] node at (2,5.5){$S_{a,c}$};
\draw[thick][fill=white] (6.5,5.5) circle[x radius=1cm, y radius=5mm] node at (8,5.5){$S_{b,c}$};
\draw[thick] (5,6.5) -- (4.25,5.8) -- (5.75,5.8) -- (5,6.5);
\end{tikzpicture}
\caption{The generalized pyramid graph $GP(r,s,t)$}\label{gra-4}
\end{figure}

They further explored the $e$-positivity of $(claw, co\text{-}diamond, F)$-free graphs where $F$ is a four-vertex graph.
The $e$-positivity of $(claw, co\text{-}diamond, F)$-free graphs is unknown for the cases $F=C_4,\, F=2K_2$ and $F=K_4$.
Hamel, Ho\`{a}ng and Tuero showed that if a peculiar graph is $(claw, co\text{-}diamond, 2K_2)$-free, then it can be characterized as a generalized pyramid $GP(r,s,t)$, as illustrated in Figure \ref{gra-4}, where $a,\,b,\,c$ are three pairwise nonadjacent vertices, the vertices of $S_{a,b}$ ($S_{a,c}$ or $S_{b,c})$ form a clique of size $r$ (resp. $s$ or $t$), and each vertex of $S_{a,b}$ ($S_{a,c}$ or $S_{b,c}$) is adjacent to every vertex of $GP(r,s,t)$ other than $c$ (resp. $b$ or $a$).
In particular, they proposed the following problem.

\begin{prob}
Whether generalized pyramids are $e$-positive?
\end{prob}
In this paper we give an affirmative answer to the above problem.

The second part of this paper is devoted to the study of the $e$-positivity of $2K_2$-free unit interval graphs. Guay-Paquet\cite{Gau-13} proved that if unit interval graphs are $e$-positive, then any claw-free incomparability graph $G$ is $e$-positive, as conjectured by Stanley. Based on Guay-Paquet's work, Foley, Ho\`{a}ng and Merkel \cite{Fol-Hoa-Mer-19} considered the $e$-positivity of $F$-free unit interval graphs, where $F$ is a four-vertex graph. It was shown that for any
four-vertex graph $F$ other than $co\text{-}diamond$, $K_4$, $4K_1$ and $2K_2$, each $F$-free unit interval graph is $e$-positive. Foley, Ho\`{a}ng and Merkel proved some special cases of $2K_2$-free unit interval graphs are $e$-positive. Based on their work, we show that any $2K_2$-free unit interval graph is $e$-positive, which provides further evidence in favor of Stanley's conjecture.

The paper is organized as follows. In Section \ref{sect-2} we prove the $e$-positivity of generalized pyramid graphs based on the monomial expansion of the corresponding chromatic symmetric function. In Section \ref{sect-3} we prove the $e$-positivity of $2K_2$-free unit interval graphs by showing that such graphs must be co-triangle free graphs or generalized bull graphs.

\section{Generalized pyramid graphs} \label{sect-2}

This section is devoted to proving the $e$-positivity of generalized pyramid graphs $GP(r,s,t)$. By using Stanley's result on the monomial expansion of the chromatic symmetric function of a graph, we first obtain the monomial expression of $X_{GP(r,s,t)}$. Then based on the transition matrix between the monomial basis and the elementary basis, we explicitly determine the coefficients in the expansion of $X_{GP(r,s,t)}$ in terms of elementary symmetric functions. Finally, we prove that all these coefficients are nonnegative.

Now let us recall some related definitions and results. Given an integer partition $\lambda$, the monomial symmetric function $m_{\lambda}$ is defined as
\begin{align*}
m_{\lambda}&=\sum_{\alpha} x^{\alpha}
\end{align*}
where $x^{\alpha}=x_1^{\alpha_1}x_2^{\alpha_2}\cdots$ and $\alpha=(\alpha_1,\,\alpha_2,\,\ldots)$ arranges over all
distinct permutations of $\lambda=(\lambda_1,\,\lambda_2,\,\ldots)$.
If $\lambda$ has $r_i$ parts equal to $i$, we also use $\langle1^{r_1}2^{r_2}\ldots\rangle$ to represent $\lambda$.
The augmented monomial symmetric function $\tilde{m}_\lambda$
is defined as
\begin{align*}
\tilde{m}_\lambda=r_1! r_2!\cdots m_\lambda.
\end{align*}
It is clear that $\{m_{\lambda} \mid \lambda\vdash n\}$ forms a basis of homogeneous symmetric functions of degree $n$, and so does $\{\tilde{m}_{\lambda} \mid \lambda\vdash n\}$.
Let $G$ be a graph  with vertex set $V$ and edge set $E$.
By using the notion of stable partitions of $G$,
Stanley \cite{Sta-95} gave a combinatorial interpretation of the coefficients in the expansion of $X_G$ in terms of $\{\tilde{m}_{\lambda}\}$.
Recall that a stable set of $G$ is a subset $S$ of $V$ such that no two vertices of $S$ are connected by an edge, and a stable partition $\pi$ of $G$ is a set partition of $V$ such that each block of $\pi$ is a stable set. The type of $\pi$ is defined to be the integer partition obtained by rearranging the block sizes of $\pi$ in decreasing order. Stanley's result can be stated as follows.

\begin{lem}\cite[Proposition 2.4]{Sta-95}\label{lem-1}
Let $G$ be a graph with $n$ vertices and $a_\lambda$ be the number of stable partitions of $G$ of type $\lambda$. Then
$$X_G=\sum_{\lambda\vdash n}a_\lambda \tilde{m}_\lambda.$$
\end{lem}

We now consider the monomial expansion of the chromatic symmetric function
of a generalized pyramid graph $GP(r,s,t)$ in Figure \ref{gra-4}.

\begin{thm}\label{pro-1}
For any nonnegative integers $r,s,t$, we have
\begin{align}\label{equ-5}
X_{GP(r,s,t)}=&\tilde{m}_{(3,1^{r+s+t})}+(rst)\tilde{m}_{(2,2,2,1^{r+s+t-3})}+(rt+rs+st+r+s+t)\tilde{m}_{(2,2,1^{r+s+t-1})}\nonumber\\[5pt]
&+(r+s+t+3)\tilde{m}_{(2,1^{r+s+t+1})}+\tilde{m}_{(1^{r+s+t+3})}.
\end{align}
\end{thm}

\pf From Figure \ref{gra-4} we see that there exists no stable set of size greater than or equal to $4$. Moreover, there exists a unique stable set of size $3$, namely $\{a,b,c\}$. A stable set of size $2$ can only be of the form $\{a,u\}$ with $u\in S_{b,c}\cup\{b,c\}$, or  $\{b,v\}$ with $v\in S_{a,c}\cup \{a,c\}$, or $\{c,w\}$ with $w\in S_{a,b}\cup \{a,b\}$.
Thus, any admissible stable partition of $GP(r,s,t)$ is of type $(3,1^{r+s+t}),\,(2,1^{r+s+t+1}),\,(2,2,1^{r+s+t-1}),\,(2,2,2,1^{r+s+t-3})$
or $(1^{r+s+t+3})$. Moreover, we have
\begin{align*}
&a_{(3,1^{r+s+t})}=1,\\[5pt]
&a_{(2,1^{r+s+t+1})}=r+s+t+3,\\[5pt]
&a_{(2,2,1^{r+s+t-1})}=rt+rs+st+r+s+t,\\[5pt]
&a_{(2,2,2,1^{r+s+t-3})}=rst,\\[5pt]
&a_{(1^{r+s+t+3})}=1.
\end{align*}
The above formulas can be proven in the same manner. As an example we prove the fourth formula. Note that a stable partition of type $(2,2,2,1^{r+s+t-3})$ is uniquely determined by the set of three stable sets of size $2$, which can only be of the form $\{\{a,u\},\{b,v\},\{c,w\}\}$ with $u\in S_{b,c},\,v\in S_{a,c},\,w\in S_{a,b}$. It is clear that $u$ has $t$ choices, $v$ has $s$ choices and $w$ has $r$ choices. Hence the fourth formula holds. This completes the proof.
\qed

Next we shall give the expansion of $X_{GP(r,s,t)}$ in terms of elementary symmetric functions. To this end, we need to use some results concerning the transition matrix between the bases $\{m_\lambda : \lambda\vdash n\}$ and $\{e_\lambda : \lambda\vdash n\}$. Let Par($n$) denote the set of all partitions of $n$. Given two partitions $\lambda=(\lambda_1,\lambda_2,\ldots)$ and $\mu=(\mu_1,\mu_2,\ldots)$ of $\text{Par}(n)$, we say that $\mu\leq\lambda$
if
$$\mu_1+\mu_2+\cdots+\mu_i\leq\lambda_1+\lambda_2+\cdots+\lambda_i\,\,\,\,\text{for all}\,i\geq1.$$
The conjugate of $\lambda=(\lambda_1,\lambda_2,\ldots)$ is defined as the partition $\lambda'=(\lambda_1',\lambda_2',\ldots)$ where $\lambda_i'=|\{j:\lambda_j\geq i\}|$. We have the following result.

\begin{lem}\cite[Chapter 7]{Sta-99}\label{lem-2}
Let $\lambda\vdash n$. If
$$e_\lambda=\sum_{\mu\vdash n}M_{\lambda\mu}m_\mu,$$
then $M_{\lambda\mu}$ is equal to the number of $(0,1)$-matrices $A=(a_{ij})_{i,j\geq1}$ satisfying $\mathrm{row}(A)=\lambda$ and $\mathrm{col}(A)=\mu$, where $\mathrm{row}(A)$ (resp., $\mathrm{col}(A)$) is the vector of row sums (resp., column sums) of $A$. Moreover, $M_{\lambda\mu}=0$ unless $\lambda\leq \mu'$, and $M_{\lambda\lambda'}=1$.
\end{lem}

Combining Theorem \ref{pro-1} and Lemma \ref{lem-2}, we obtain the following result.

\begin{thm}\label{thm-1}
For any nonnegative integers $r,s,t$, we have
\begin{align}\label{equ-6}
X_{GP(r,s,t)}=&A\cdot e_{(r+s+t+1,1,1)}+B\cdot e_{(r+s+t,3)}+C\cdot e_{(r+s+t+1,2)}\nonumber\\
&\quad +D\cdot e_{(r+s+t+2,1)}+E\cdot e_{(r+s+t+3)},
\end{align}
where
\begin{align*}%\label{equ-7}
A=&(r+s+t)!,\nonumber\\[5pt]
B=&(r+s+t-3)!\cdot 6rst,\nonumber\\[5pt]
C=&(r+s+t-3)!\cdot 2(r + s + t - 1)\nonumber\\
&\cdot[(r^2 s + r s^2-2 r s) + (r t^2+ r^2 t - 2 r t)+(s^2 t + s t^2- 2 s t)],\nonumber\\[5pt]
D=&( r + s + t - 2)!\cdot [( r^4+r^3-2 r^2 )+(3 r^2 s- 2 r s)+(3 r s^2-2 s^2)\nonumber\\
&+(3 r^2 t- 2 r t) + (9 r s t- 2 s t )+ (3 r t^2- 2 t^2 )+ 3 s^2 t +5 r s^2 t \nonumber\\
& + 2 s^3 t + 5 r^2 s t + 2 r^3 t  + 2 r^2 t^2 + 3 s t^2 +5 r s t^2 + 2 s^2 t^2 \nonumber\\
&+ t^3 + 2 r t^3 + 2 s t^3 + t^4 + 2 r^3 s+ 2 r^2 s^2 + s^3 + 2 r s^3+ s^4 ],\nonumber\\[5pt]
E=&( r + s + t - 1)!\cdot (3 + r + s + t)(r + s)(r + t)(s + t).
\end{align*}
\end{thm}

\pf
Let $i=r+s+t$ and $P=\{(2^3,1^{i-3}),\,(3,1^{i}),\,(2^2,1^{i-1}),\,
(2,1^{i+1}),\,(1^{i+3})\}$. In order to give the elementary expansion of $X_{GP(r,s,t)}$, by Theorem \ref{pro-1} and Lemma \ref{lem-2} it suffices to consider the monomial expansion of those $e_{\lambda}$'s such that $\lambda'\leq \mu$ for some $\mu\in P$. It is straightforward to verify that the set of such partitions $\lambda$ is composed of $\{(i,3),\,(i+1,1,1),\,(i+1,2),\,(i+2,1),\,(i+3)\}$.
By Lemma \ref{lem-2}, we get
\begin{align}
&e_{(i,3)}=m_{(2,2,2,1^{i-3})}+(i-1)m_{(2,2,1^{i-1})}+
\binom{i+1}{2}m_{(2,1^{i+1})}+\binom{i+3}{3}m_{(1^{i+3})},\label{equ-1}\\[5pt]
&e_{(i+1,1,1)}=m_{(3,1^i)}+(2i+3)m_{(2,1^{i+1})}+2 m_{(2,2,1^{i-1})}+2\binom{i+3}{2}m_{(1^{i+3})},\label{equ-2}\\[5pt]
&e_{(i+1,2)}=m_{(2,2,1^{i-1})}+(i+1)m_{(2,1^{i+1})}+\binom{i+3}{2}m_{(1^{i+3})},\label{equ-3}\\[5pt]
&e_{(i+2,1)}=m_{(2,1^{i+1})}+(i+3)m_{(1^{i+3})},\label{equ-4}\\[5pt]
&e_{i+3}=m_{(1^{i+3})}.\label{equ-0}
\end{align}
The above formulas are easy to prove. As an example we prove that the coefficient of $m_{(2,1^{i+1})}$ in $e_{(i+1,2)}$ is $(i+1)$. By Lemma \ref{lem-2}, we only need to count the number of $(0,1)$-matrices $A=(a_{pq})_{p,q\geq1}$ with $\text{row}(A)=(i+1,2)$ and $\text{col}(A)=(2,1^{i+1})$. Since $\text{row}(A)=(i+1,2)$, there are $i+1$ entries equal to $1$ in the first row of matrix $A$ and there are $2$ entries equal to $1$ in the second row. Since $\text{col}(A)=(2,1^{i+1})$, we must have $a_{11}=a_{21}=1$ and $a_{pq}=0$ for $p \geq 3$ or $q\geq i+3$. Moreover, the submatrix
$$
\begin{pmatrix}
a_{12} & a_{13} & \cdots & a_{1,i+2}\\
a_{22} & a_{23} & \cdots & a_{2,i+2}\\
\end{pmatrix}
$$
can be any $2\times (i+1)$ matrix composed of $i$ column vectors $\left({1\atop 0}\right)$'s and one column vector $\left({0\atop 1}\right)$.
Hence we have $M_{(i+1,2),(2,1^{i+1})}=i+1$.

By using the above $m$-expansion formulas we can get the $e$-expansion of those monomial symmetric functions appearing in \eqref{equ-5}. Substituting the resulted $e$-expansion formulas into \eqref{equ-5}, we complete the proof.
\qed

We proceed to prove the main result of this section.

\begin{thm}\label{thm-gen-pyr}
For any nonnegative integers $r,s,t\geq 0$ the generalized pyramid graph $GP(r,s,t)$ is $e$-positive.
\end{thm}

\pf Note that if $r=s=t=0$, then $X_{GP(r,s,t)}=e_1^3$, which is obviously $e$-positive. If only two of $r,s,t$ are zero, then $GP(r,s,t)$ belongs to the class of $e$-positive graphs studied by Hamel, Ho\`{a}ng and Tuero, see \cite[Lemma 9]{Ham-Hoa-Tue-19}. If only one of $r,s,t$ is zero, then $GP(r,s,t)$ is a generalized bull graph, whose positivity is already known, see Foley, Ho\`{a}ng and Merkel \cite[Theorem 11]{Fol-Hoa-Mer-19} and Cho, Huh \cite[Theorem 3.6]{Cho-Huh-18}.

From now on we assume that $r,s,t$ are positive integers. In order to show the $e$-positivity of $X_{GP(r,s,t)}$, it suffices to show that the coefficients $A,B,C,D,E$ in \eqref{equ-6} are nonnegative by Theorem \ref{thm-1}. Clearly, $A$, $B$ and $E$ are always nonnegative.

We continue to prove $C\geq0$.
Since $r,s\geq1$, we have
\[r^2s+rs^2-2rs\geq r^2+s^2-2rs\geq 0,\]
Similarly, we have
\[r^2t+rt^2-2rt\geq 0,\]
and
\[st^2+s^2t-2st\geq 0.\]
Therefore, $C\geq 0$.

Finally, we prove that $D\geq 0$. Since $r,s,t\geq 1$, it is straightforward to verify that
$r^4+r^3 -2 r^2, 3 r^2 s- 2 r s, 3 r s^2-2 s^2, 3 r^2 t- 2 r t, 9 r s t- 2 s t, 3 r t^2- 2 t^2$ are all nonnegative.
Thus, $D\geq 0$. This completes the proof.\qed

\section{\texorpdfstring{$2K_2$},-free unit interval graphs} \label{sect-3}

The aim of this section is to prove that $2K_2$-free unit interval graphs are $e$-positive.
Our proof is based on the characterization of $2K_2$-free unit interval graphs due to Hempel and Kratsch \cite{Hem-Kra-02}, who
actually gave a characterization of a larger family of graphs. Using their result, we show that $2K_2$-free unit interval graphs can only be either co-triangle-free graphs or generalized bull graphs, which are already known to be $e$-positive.

Let us first recall some related definitions and results. A co-triangle means a stable set of size $3$. Stanley and Stembridge \cite{Sta-Ste-93}
proved the $e$-positivity of the complement of bipartite graphs, which are a special class of co-triangle-free graphs.
Stanley \cite{Sta-95} gave a different proof of their result, and his arguments can also be applied to the following general case.

\begin{lem}\cite[Exercise 7.47]{Sta-99}\label{lem-stanley-1}
If $G$ is a co-triangle-free graph, then $X_G$ is $e$-positive.
\end{lem}

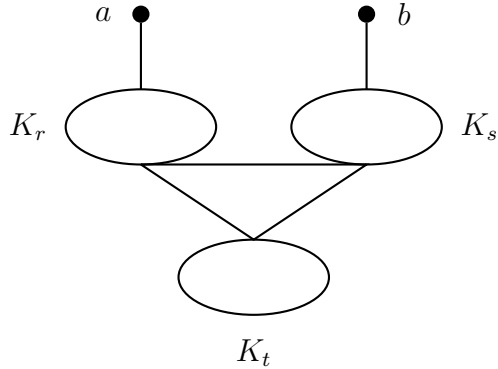
\begin{figure}[ht]
\center
\begin{tikzpicture}
[place/.style={thick,fill=black!100,circle,inner sep=0pt,minimum size=2mm,draw=black!100}]
\node [place] (v2) at (7,6.5) {};
\node [place] (v1) at (4,6.5) {};
\draw [thick] (4,5) ellipse (1 and 0.5);
\draw [thick] (7,5) ellipse (1 and 0.5);
\draw [thick] (5.5,3) ellipse (1 and 0.5);
\draw[thick] (v1) -- (4,5.5);
\draw[thick] (v2) -- (7,5.5);
\draw[thick] (4,4.5) -- (7,4.5) -- (5.5,3.5) -- (4,4.5);
\node at (8.5,5) {$K_s$};
\node at (2.5,5) {$K_r$};
\node at (5.5,2) {$K_t$};
\node at (3.5,6.5) {$a$};
\node at (7.5,6.5) {$b$};
\end{tikzpicture}
\caption{The generalized bull graph $GB(r,s,t)$}\label{gra-2}
\end{figure}

The generalized bull graphs were introduced by Foley, Ho\`{a}ng and Merkel \cite{Fol-Hoa-Mer-19}, but their $e$-positivity was first studied by Cho and Huh \cite{Cho-Huh-18}. A generalized bull graph can be characterized as in Figure \ref{gra-2}, where $K_r,\,K_s,\,K_t$ form a clique of size $r+s+t$, $a$ is adjacent to each vertex of $K_r$, and $b$ is adjacent to each vertex of $K_s$. We denote such a graph by $GB(r,s,t)$.

Cho and Huh \cite{Cho-Huh-18} obtained the following result.

\begin{lem}\cite[Theorem 3.6]{Cho-Huh-18}\label{lem-gen-bull}
For any positive integers $r,s,t$, the generalized bull graph ${GB(r,s,t)}$ is $e$-positive.
\end{lem}

Note that Cho and Huh proved the above result based on the Schur expansion of $X_{GB(r,s,t)}$.
To be self-contained, we would like to give a new proof, which parallels that of Theorem \ref{thm-gen-pyr}.

\noindent\textit{Proof of Lemma \ref{lem-gen-bull}.}
We first give the monomial expansion of $X_{GB(r,s,t)}$.
Using the same method as in the proof of Theorem \ref{pro-1}, we get that
\begin{align}\label{equ-8}
X_{GB(r,s,t)}=&t\cdot\tilde{m}_{(3,1^{r+s+t-1})}+(t(t-1)+tr+sr+st)\cdot\tilde{m}_{(2,2,1^{r+s+t-2})}\nonumber\\[5pt]
&+(1+2t+s+r)\cdot\tilde{m}_{(2,1^{r+s+t})}+\tilde{m}_{(1^{r+s+t+2})}.
\end{align}

Setting $k=r+s+t$ and $i=k-1$ in \eqref{equ-2}, \eqref{equ-3}, \eqref{equ-4} and \eqref{equ-0}, and then substituting these four equations into \eqref{equ-8}, we obtain
\begin{align*}
X_{GB(r,s,t)}=&(r+s+t-2)!\cdot[(r + s + t - 1)t\cdot e_{(r+s+t,1,1)}+2rs\cdot e_{(r+s+t,2)}\\[5pt]
&+(r^3 + r^2 s + r s^2 + s^3 + 2 r^2 t + 2 r s t + 2 s^2 t + r t^2 + s t^2
- r - s)\cdot e_{(r+s+t+1,1)}\\[5pt]
&+(r + s + t + 2)(r + s + t - 1) r s\cdot e_{(r + s + t + 2)}].
\end{align*}

Since $r,s,t\geq 1$, the $e$-positivity of $X_{GB(r,s,t)}$ is obvious. \qed

We proceed to recall Hempel and Kratsch's characterization of $2K_2$-free unit interval graphs. As will be shown below,
$2K_2$-free unit interval graphs are a special class of (claw, AT)-free graphs. Recall that an interval graph is formed from a set of intervals on the real line, with a vertex for each interval and an edge between vertices whose intervals intersect. A unit interval graph is an interval graph for which each of its intervals has unit length. It is well known that unit interval graphs must be claw-free and $C_4$-free. The notion of AT-free graphs was introduced by Lekkerkerker and Boland \cite{Lek-Bol-62}. A co-triangle in a graph $G$ is called an asteroidal triple, denoted by AT for short, if for any pair of its vertices there exists a path between them which non-intersects with the neighborhood of the third vertex. It has been shown in \cite{Lek-Bol-62} that interval graphs are exactly the class of chordal AT-free graphs, where a graph is a chordal graph if every induced cycle in the graph have exactly three vertices. Meanwhile, unit interval graphs have been shown to be exactly the class of claw-free interval graphs \cite{Gar-07}. Hence, $2K_2$-free unit interval graphs are equivalent to ($2K_2$, claw, AT)-free chordal graphs. Given a graph $G$ with vertex set $V$ and edge set $E$ and a pair of vertices $u$ and $v$, let $\alpha(G)$ denote the maximum size of stable sets and let $d(u,v)$ denote the number of edges of the shortest path between $u$ and $v$. For any vertex $w\in V$, let $N_i(w)=\{x\in V\mid d(x,w)=i\}$ and $[N_i(w)]$ denote the induced subgraph on $N_i(w)$. In particular, $N_1(w)$ is the neighborhood of $w$, denoted by $N(w)$. With these notations,
Hempel and Kratsch's characterization of (claw, AT)-free graphs can be stated as follows.

\begin{lem}\cite[Lemma 6]{Hem-Kra-02}\label{lem-AT}
For any connected (claw, AT)-free graph $G$, there exists a vertex $w$ such that $\alpha([N(w)])\leq 2$ and for any $i\geq 2$ each $[N_i(w)]$ is a clique (which might be empty).
\end{lem}

It is well known that $X_{G\uplus H}=X_G X_H$, where $G\uplus H$ is a disjoint union of graphs $G$ and $H$.
Given a $2K_2$-free unit interval graph $G$, it is clear that every connected component of $G$ is also a $2K_2$-free unit interval graph. Thus in order to study the $e$-positivity of $X_G$, we may assume that $G$ is connected.
Based on the above result, we could give a characterization of connected $2K_2$-free unit interval graphs.

\begin{coro}\label{coro-2k2}
If $G$ is a connected $2K_2$-free unit interval graph, then there exists a vertex $w$ such that $\alpha([N(w)])\leq 2$, $|N_2(w)|$ is a clique, $|N_3(w)|\leq 1$, and $N_i(w)=\emptyset$ for any $i\geq 4$. Moreover, if $[N(w)]$ is connected, $|N_3(w)|=0$ and $\alpha([N(w)])=2$, then $|N_2(w)| \leq 2$ and $[N(p)\cap N(w)]$ is a clique for any $p\in N_2(w)$.
\end{coro}

\noindent \textit{Proof of Corollary \ref{coro-2k2}.} Since $G$ is a $2K_2$-free unit interval graph, thus it must be (claw, AT)-free, as mentioned before Lemma \ref{lem-AT}. Thus, there exists $w$ such that $\alpha([N(w)])\leq 2$ and for any $i\geq 2$ each $[N_i(w)]$ is a clique.

We proceed to show that $|N_3(w)|\leq 1$ and $N_i(w)=\emptyset$ for any $i\geq 4$. We first show that $N_i(w)=\emptyset$ for any $i\geq 4$. Otherwise, if $N_i(w)\neq\emptyset$ for some $i\geq 4$, then $N_{j}(w)\neq\emptyset$ for any $1\leq j\leq i-1$. Thus there exist $x\in N(w)$, $y\in N_{i-1}(w)$ and $z\in N_i(w)$ such that the set $\{w,\,x,\,y,\,z\}$ induces a $2K_2$, a contradiction.
We next show that $|N_3(w)|\leq 1$. Otherwise if $|N_3(w)|>1$, then there exist $u,v\in N_3(w)$ such that $uv\in E$, since $[N_3(w)]$ is a clique. Then for any $x$ in $N(w)$, the set $\{w,\,x,\,u,\,v\}$ induces a $2K_2$, a contradiction. Hence $|N_3(w)|\leq 1$.

It remains to show that if $[N(w)]$ is connected, $|N_3(w)|=0$ and $\alpha([N(w)])=2$, then $|N_2(w)| \leq 2$ and $[N(p)\cap N(w)]$ is a clique for any $p\in N_2(w)$. Note that a $2K_2$-free unit interval graph must be $C_4$-free.
We first show that $[N(p)\cap N(w)]$ is a clique for any $p\in N_2(w)$.  Suppose to the contrary
there exist some $p\in N_2(w)$ and non-adjacent $a,b\in N(p)\cap N(w)$.
Thus $\{p,\,a,\,b,\,w\}$ induces a $C_4$, a contradiction. We next show that $|N_2(w)| \leq 2$. Suppose that $|N_2(w)|=s$.
We claim that for any $a\in N(w)$ there are at least $s-1$ vertices in $N_2(w)$ which are adjacent to $a$, namely, $|N(a)\cap N_2(w)|\geq s-1$. Suppose to the contrary there exist some $a\in N(w)$ and $x,y\in N_2(w)$ such that neither $x$ nor $y$ is adjacent to $a$, and thus $\{x,\,y,\,a,\,w\}$ induces a $2K_2$ in $G$ since $[N_2(w)]$ is a clique, a contradiction. Since $\alpha([N(w)])=2$, there exist $a,b\in N(w)$ such that $a$ and $b$ are not adjacent.
Moreover, $a,b$ can not be adjacent to the same vertex $x$ in $N_2(w)$ for otherwise the set $\{x,\,a,\,b,\,w\}$ induces a $C_4$, a contradiction.
 This means that
$$(N(a)\cap N_2(w))\cap(N(b)\cap N_2(w))=\emptyset.$$
Hence
\[s=|N_2(w)|\geq |N(a)\cap N_2(w)|+|N(b)\cap N_2(w)|\geq (s-1)+(s-1),\]
yielding $s\leq 2$. Hence $|N_2(w)| \leq 2$. This completes the proof.
\qed

We would like to point out that the first part of Corollary \ref{coro-2k2} was already known to Foley, Ho\`{a}ng and Merkel \cite{Fol-Hoa-Mer-19}, and the second part tells more information of a $2K_2$-free unit interval graph $G$. In fact,
if more constraints are added, we could get a clearer characterization of $G$. The following result will be used to check the $e$-positivity of some special $2K_2$-free unit interval graphs.

\begin{coro}\label{coro-2k2-1}
 Given a connected $2K_2$-free unit interval graph $G$, let $w$ be given as in Corollary \ref{coro-2k2}.
 Suppose that $[N(w)]$ is connected, $|N_2(w)|=1$, $|N_3(w)|=0$ and $\alpha([N(w)])=2$.
 If we let $N_2(w)=\{p\}$, $A=N(p)\cap N(w)$ and $B=N(w)\setminus A$, then $|N(a)\cap B|\geq |B|-1$ and $[N(a)\cap B]$ is a clique for any $a\in A$.
\end{coro}

\pf Let first prove that $|N(a)\cap B|\geq |B|-1$ for any $a\in A$.
Suppose the contrary. There exist $a\in A$ and $b_1,\,b_2\in B$ such that $b_1$ and $b_2$ are not adjacent to $a$.
If $b_1$ and $b_2$ are not adjacent in $G$, then $\{a,\,b_1,\,b_2\}$ is a stable set, contradicting $\alpha([N(w)])=2$.
If $b_1$ and $b_2$ are adjacent, then $\{a,\,p,\,b_1,\,b_2\}$ induces a $2K_2$, a contradiction. Thus $a$ is adjacent to at least $|B|-1$ vertices in $B$.
Next we show that  $[N(a)\cap B]$ is a clique for any $a\in A$. Suppose to the contrary there exists some $a\in A$ and non-adjacent $b,b'\in N(a)\cap B$. Note that the set $\{a,\,p,\,b,\,b'\}$ induces a claw, a contradiction. This completes the proof. \qed

Finally we come to the main result of this section.

\begin{thm}
If $G$ is a $2K_2$-free unit interval graph, then $X_G$ is $e$-positive.
\end{thm}

\pf Without loss of generality, we may assume that $G$ is connected. By Corollary \ref{coro-2k2}, there are thus six cases to check:
\begin{itemize}
\item[(1)] $[N(w)]$ is not connected;
\item[(2)] $[N(w)]$ is connected and $|N_3(w)|=1$;
\item[(3)] $[N(w)]$ is connected, $|N_3(w)|=0$ and $\alpha([N(w)])=1$;
\item[(4)] $[N(w)]$ is connected, $|N_3(w)|=0$, $\alpha([N(w)])=2$ and $|N_2(w)|=2$;
\item[(5)] $[N(w)]$ is connected, $|N_3(w)|=0$, $\alpha([N(w)])=2$ and $|N_2(w)|=1$;
\item[(6)] $[N(w)]$ is connected, $|N_3(w)|=0$, $\alpha([N(w)])=2$ and $|N_2(w)|=0$;
\end{itemize}
where $w$ is given as in Corollary \ref{coro-2k2}.

Foley, Ho\`{a}ng and Merkel \cite{Fol-Hoa-Mer-19} showed that the theorem is true for the first three cases, for which they showed that $G$ must be a co-triangle free graph or a generalized bull graph. We only need to consider the remaining three cases.

Let us first consider Case (6). In this case, it is clear that $G$ is co-triangle-free. Thus $X_G$ is $e$-positive by Lemma \ref{lem-stanley-1}.

\begin{figure}[ht]
\center
\begin{tikzpicture}
[place/.style={thick,fill=black!100,circle,inner sep=0pt,minimum size=2mm,draw=black!100}]
\node [place,label=right:$q$] (v2) at (7,6.5) {};
\node [place,label=left:$p$] (v1) at (4,6.5) {};
\node [place,label=below:$w$] (v3) at (5.5,3.5) {};
\draw [thick] (4,5) ellipse (1 and 0.5) ;
\node at (2.75,5) {$A$};
\node at (8.25,5) {$B$};
\draw [thick] (7,5) ellipse (1 and 0.5);
\draw [thick] (4,4.5) -- (v3) -- (7,4.5);
\draw [thick] (4,5.5) -- (v1) -- (v2) -- (7,5.5);
\node (v4) at (4.85,4.9) {};
\node (v5) at (6.15,5) {};
\draw[thick] [densely dashed](v2) -- (v4) ;
\node (v6) at (4.9,5) {};
\draw[thick] [densely dashed](v6) -- (v5) ;
\end{tikzpicture}
\caption{A co-triangle-free graph.}\label{gra-5}
\end{figure}
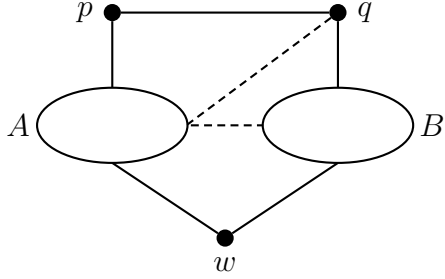

Next we consider Case (4). Set $N_2(w)=\{p,q\}$, $A=N(p)\cap N(w)$ and $B=N(w)\setminus A$. By Corollary \ref{coro-2k2}, $[A]$ is a clique. We claim that any vertex $b\in B$ is adjacent to $q$. Otherwise if there exists some $b\in B$ such that $q$ and $b$ are not adjacent, then $\{p,\,q,\,b,\,w\}$ induces a $2K_2$, a contradiction. Hence all vertices of $B$ are adjacent to $q$. By Corollary \ref{coro-2k2} the induced subgraph $[N(q)\cap N(w)]$ is a clique and hence $[B]$ is a clique. Thus $G$ can be characterized as a co-triangle-free graph, as depicted in Figure \ref{gra-5}, where the dashed lines represent that there may exist some edges between $A$ and $B$, as well as between $q$ and $A$. Again by Lemma \ref{lem-stanley-1}, we obtain the $e$-positivity of $X_G$.

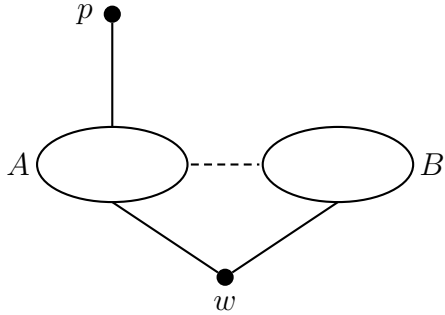
\begin{figure}[ht]
\center
\begin{tikzpicture}
[place/.style={thick,fill=black!100,circle,inner sep=0pt,minimum size=2mm,draw=black!100}]
\node [place,label=left:$p$] (v1) at (4,7) {};
\node [place,label=below:$w$] (v3) at (5.5,3.5) {};
\draw [thick] (4,5) ellipse (1 and 0.5) ;
\node at (2.75,5) {$A$};
\node at (8.25,5) {$B$};
\draw [thick] (7,5) ellipse (1 and 0.5);
\draw [thick] (4,4.5) -- (v3) -- (7,4.5);
\draw [thick] (v1) -- (4,5.5);
\node (v2) at (4.9,5) {};
\node (v4) at (6.1,5) {};
\draw [thick][densely dashed] (v2) edge (v4);
\end{tikzpicture}
\caption{A co-triangle-free graph.}\label{gra-6}
\end{figure}

Finally, we prove that the theorem holds for Case (5). Now set $N_2(w)=\{p\}$, $A=N(p)\cap N(w)$ and $B=N(w)\setminus A$.
By Corollary \ref{coro-2k2}, $[A]$ is a clique. If $[B]$ is a clique, then it is easy to see that $G$ is co-triangle-free, see Figure \ref{gra-6}, where the dashed lines represent that there may exist some edges between $A$ and $B$. Hence $X_G$ is $e$-positive by Lemma \ref{lem-stanley-1}.

\begin{figure}[ht]
\center
\begin{tikzpicture}
[place/.style={thick,fill=black!100,circle,inner sep=0pt,minimum size=2mm,draw=black!100}]
\node [place,label=above:$p$] (v1) at (5,9.5) {};
\node [place,label=left:$x$] (v4) at (2.5,5.5) {};
\node [place,label=right:$y$] (v2) at (7.5,5.5) {};
\node [place,label=below:$w$] (v3) at (5,3) {};
\draw [thick](3.5,7.5) ellipse (1 and 0.5) ;
\draw [thick](6.5,7.5) ellipse (1 and 0.5);
\draw [thick](5,5.5) ellipse (1 and 0.5);
\node at (2.25,7.5) {$A_1$};
\node at (7.75,7.5) {$A_2$};
\node at (5.25,4.75) {$A_3$};
\draw [thick](3.5,8) -- (v1) -- (6.5,8);
\draw [thick](4.5,7.5) -- (5.5,7.5);
\draw [thick](v2) -- (v3) -- (5,5) -- (v3) -- (v4);
\draw [thick]plot[smooth, tension=.7] coordinates {(v3) (2,5) (2.75,7.15)};
\draw [thick]plot[smooth, tension=.7] coordinates {(5,3) (8,5) (7.25,7.15)};
\draw [thick](v4) -- (3.5,7) -- (4.51,5.92);
\draw [thick](v2) -- (6.5,7) -- (5.51,5.92);
\draw [thick](v4) -- (4,5.5);
\draw [thick](v2) -- (6,5.5);
\end{tikzpicture}
\caption{{x,y,p} induces an asteroidal triple.}\label{gra-7}
\end{figure}
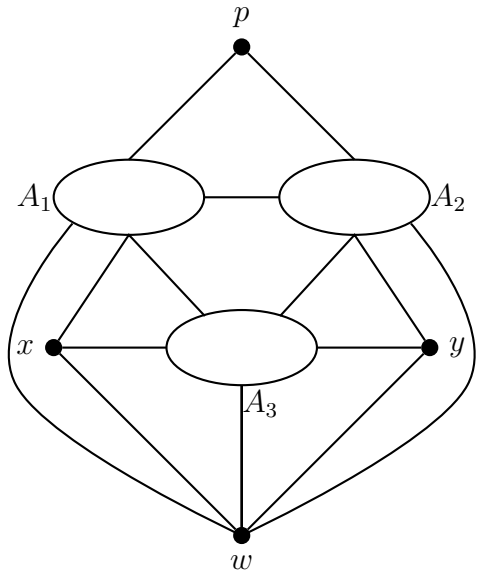

From now on we assume that $[B]$ is not a clique. Thus there exist non-adjacent vertices $x,y\in B$.
Now set $A_1= N(x)\cap A$, $A_2=N(y)\cap A$ and $A_3=B\setminus\{x\cup y\}$. We claim that
either $A_1 =\emptyset$ or $A_2=\emptyset$. Suppose to the contrary that $A_1\neq \emptyset$ and $A_2\neq \emptyset$.
Now we have $A_1\cap A_2=\emptyset$, otherwise there exists $a\in A_1\cap A_2$ such that $\{a,\,x,\,y,\,p\}$ induces a claw, a contradiction. Moreover, we have $A_1 \cup A_2=A$, otherwise there exists $b\in A\setminus(A_1\cup A_2)$ such that $\{b,\,x,\,y\}$ is a stable set, contradicting $\alpha([N(w)])\leq 2$. By Corollary \ref{coro-2k2-1}, both $[\{x\}\cup A_1\cup A_3]$ and $[\{y\}\cup A_2 \cup A_3]$ are cliques. A little thought shows that $\{x,\,y,\,p\}$ is an asteroidal triple, as shown in Figure \ref{gra-7}. While this contradicts the fact that $G$ is AT-free. Thus at least one of $A_1$ and $A_2$ is empty.

Without loss of generality, we may assume that $A_1=\emptyset$.
We proceed to show that $N(w)\setminus\{x\}$ induces a clique.
By noting that $N(w)=A\cup B$ and $A$ is a clique, it suffices to show that each $a\in A$ and each $z\in B\setminus \{x\}$ are adjacent and $[B\setminus\{x\}]$ is a clique.
For the former assertion, assume to the contrary that there exist non-adjacent $a\in A$ and $z\in B\setminus \{x\}$. If $x,z$ are not adjacent, then $\{w,\,x,\,a,\,z\}$ induces a claw, a contradiction. If $x,z$ are adjacent, then $\{a,\,p,\,x,\,z\}$ induces a $2K_2$, again a contradiction. Hence $a$ and $z$ are adjacent. For the latter assertion, assume to the contrary that  there exist non-adjacent vertices $b_1,b_2\in B\setminus\{x\}$ such that for any $a\in A$ the set $\{a,\,p,\,b_1,\,b_2\}$ induces a claw, a contradiction. Thus $[N(w)\setminus\{x\}]$ is a clique.

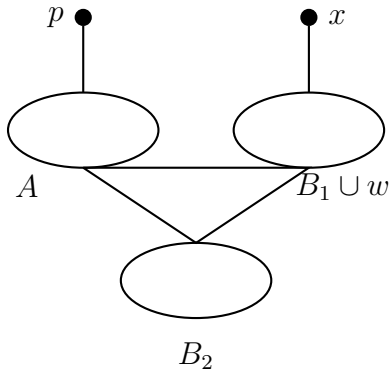
\begin{figure}[ht]
\center
\begin{tikzpicture}
[place/.style={thick,fill=black!100,circle,inner sep=0pt,minimum size=2mm,draw=black!100}]
\node [place,label=right:$x$] (v2) at (7,6.5) {};
\node [place,label=left:$p$] (v1) at (4,6.5) {};
\draw [thick] (4,5) ellipse (1 and 0.5);
\node at (3.25,4.25) {$A$};
\draw [thick] (7,5) ellipse (1 and 0.5);
\node at (7.45,4.25) {$B_1\cup w$};
\draw [thick] (5.5,3) ellipse (1 and 0.5);
\node at (5.5,2) {$B_2$};
\draw[thick] (v1) -- (4,5.5);
\draw[thick] (v2) -- (7,5.5);
\draw[thick] (4,4.5) -- (7,4.5) -- (5.5,3.5) -- (4,4.5);
\end{tikzpicture}
\caption{A generalized bull graph.}\label{gra-8}
\end{figure}

If we set $B_1=N(x)\cap (B\setminus\{x\})$ and $B_2=B\setminus\{x\cup B_1\}$, then $G$ can be considered as a generalized bull graph, see Figure \ref{gra-8}.
Thus for case (5) if $[B]$ is not a clique, the graph $G$ is also $e$-positive by Lemma \ref{lem-gen-bull}.

Combining the above all cases, we complete the proof. \qed

\noindent \textbf{Acknowledgments.}
This work is supported in part by the Fundamental Research Funds for the Central Universities and the National Science Foundation of China (Nos. 11522110, 11971249). We thank Ethan Y.H. Li for stimulating discussions.

\end{document}